\documentclass[11pt]{amsart}
\usepackage{dsfont}
\usepackage{times}
\usepackage{enumerate}
\usepackage[euler-digits]{eulervm}
\theoremstyle{plain}
\usepackage{graphicx,color} 
\usepackage{mathtools}

\usepackage[hidelinks]{hyperref}
\usepackage{subcaption}
\usepackage{amsmath, amssymb, graphics}
\usepackage{float}

\newtheorem{theorem}{Theorem}[section] 
\newtheorem*{theorem*}{Theorem}
\newtheorem{cor}{Corollary}[section]

\newtheorem{prop}{Proposition}[section]

\newtheorem{remark}{Remark}[section]

\newtheorem*{conjecture*}{Conjecture}

\begin{document}
	
\title{Hyperbolic Geometry and the Helfrich Functional}

\author{Bennett Palmer and \'Alvaro P\'ampano}

\date{\today}

\maketitle

\begin{abstract}{The Helfrich model is a fundamental tool for determining the morphology of biological membranes. We relate the geometry of an important class of its equilibria to the geometry of sessile and pendant drops in the hyperbolic space ${\bf H}^3$. When the membrane surface meets the ideal boundary of hyperbolic space, a modification of the regularized area functional is related to the construction of closed equilibria for the Helfrich functional in ${\bf R}^3$.}\\
	
\noindent K{\tiny EY} W{\tiny ORDS}.\, Helfrich Energy, Hyperbolic Space, Sessile and Pendant Drops, Reduced Membrane Equation, Renormalized Area.
	
\noindent MSC C{\tiny LASSIFICATION} (2020).\, 49Q10, 53A05, 58E50.
\end{abstract}

\section{Introduction}

The Helfrich functional
\begin{equation*} 
	\mathcal{H}_{a,c_o,b}[\Sigma]:=\int_\Sigma \left(a\left[H+c_o\right]^2+bK\right)d\Sigma\:,
\end{equation*}
is the principal mechanism for determining the morphology of lipid bilayers \cite{H}. In this model, bilipid membranes are modeled by mathematical surfaces $\Sigma$ with mean curvature $H$ and Gaussian curvature $K$. The parameters $a>0$, $c_o, b\in{\bf R}$ of the functional $\mathcal{H}_{a,c_o,b}$ are physical constants depending on the composition of the membrane. The Helfrich functional $\mathcal{H}_{a,c_o,b}$ is a type of elastic energy, a particular case of which, $\mathcal{H}_{a,0,-a}$, is well known to be conformally invariant \cite{B}. In this case it is not surprising that, when the surface $\Sigma$ is contained in a half-space, hyperbolic geometry plays a role in determining the values of the functional.

In the current paper, we study relationships between equilibria for $\mathcal{H}_{a,c_o,b}$ and hyperbolic geometry for general values of $c_o\in{\bf R}$. In \cite{PP2}, the authors showed that there is a special class of equilibria for the Helfrich functional which, after certain normalizations,  satisfy the reduced membrane equation
\begin{equation*}
H+c_o=-\frac{\nu_3}{z}\:.
\end{equation*}
In this equation, $z$ is a coordinate of the surface (which we take to be the vertical one) and $\nu_3$ is the vertical component of the unit normal to the surface. Assuming that the surface lies in the upper half-space, if the background metric is changed from the Euclidean to the hyperbolic one, the mean curvature of $\Sigma$ in ${\bf H}^3$ is ${\widehat H}=Hz+\nu_3$, so the reduced membrane equation becomes ${\widehat H}=-c_oz$. Surfaces in Euclidean space with mean curvature depending linearly on their height have been used to model liquid drops in a vertical gravitational field and, hence, they have been widely studied as capillary surfaces with gravity \cite{F,W}. Similarly, surfaces satisfying the reduced membrane equation are equilibria for the functional
\begin{equation*}
	\mathcal{G}[\Sigma]:={\widehat{ \mathcal A}}[\Sigma]-2c_o{\widehat {\mathcal U}}[\Sigma]\:,
\end{equation*}
where ${\widehat{ \mathcal A}}$ denotes the hyperbolic area and ${\widehat {\mathcal U}}$ the (hyperbolic) gravitational potential energy of the volume enclosed by the surface. When such a surface extends to the ideal boundary of the hyperbolic space ${\bf H}^3$ in such a way that it intersects this boundary in a right angle, then the functional $\mathcal{G}$ may become infinite, however it possesses a regularization $\mathcal{G}_R$ for which the inequality
\begin{equation*}
-\mathcal{G}_R[\Sigma]\le \mathcal{H}_{1,c_o,0}[\Sigma]
\end{equation*}
holds with equality holding for equilibria of $\mathcal{G}_R$.

\section{Preliminaries}

Let $(x,y,z)$ denote the canonical coordinates of the Euclidean $3$-space ${\bf R}^3$ and let $E_i$, $i=1,2,3$, be the constant unit vector fields in the direction of the coordinate axes. Consider a connected and oriented compact surface $\Sigma$ with boundary and assume that $\Sigma$ is smoothly embedded in ${\bf R}^3$. Denote by $\nu:\Sigma\longrightarrow{\bf S}^2$ the unit normal vector field defined by the property that it points out of any convex domain.

The \emph{Helfrich energy} \cite{H} is the functional
\begin{equation}\label{H}
	\mathcal{H}_{a,c_o,b}[\Sigma]:=\int_\Sigma\left(a\left[H+c_o\right]^2+b K\right)d\Sigma\,,
\end{equation}
where $H$ and $K$ denote the mean and Gaussian curvatures of $\Sigma$, respectively. The parameters $a>0$, $c_o$, and $b$ are real constants. The specific case $a=-b$ and $c_o=0$ gives rise to the classical Willmore energy \cite{Wi}. The functional $\mathcal{H}_{a,0,-a}$ is conformally invariant \cite{B}, a property that no longer holds for other choices of the energy parameters.

The Euler-Lagrange equation associated to the Helfrich energy \eqref{H} is the fourth order nonlinear elliptic partial differential equation
\begin{equation}\label{ELH}
	\Delta H+2(H+c_o)\left(H\left[H-c_o\right]-K\right)=0\,,
\end{equation}
where $\Delta$ denotes the Laplace-Beltrami operator. In \cite{PP2}, the authors showed that (after a suitable rigid motion and translation of the vertical coordinate) a special class of solutions to \eqref{ELH} is given by the second order equation
\begin{equation}\label{RME}
	H+c_o=-\frac{\nu_3}{z}\,,
\end{equation}
where $\nu_3:=\nu\cdot E_3$ is the vertical component of the unit normal. The equation \eqref{RME} is known as the \emph{reduced membrane equation}.

For any part of the surface $\Sigma$ which lies in the upper\footnote{After the transformation $z\to-z$, which preserves \eqref{RME}, the same is true for the part of $\Sigma$ lying in the lower half-space.} half-space ${\bf R}^3_+=\{z>0\}$, the reduced membrane equation \eqref{RME} can be interpreted as the Euler-Lagrange equation for the functional
\begin{equation}\label{G}
	\mathcal{G}[\Sigma]:=\int_\Sigma \frac{1}{z^2}\,d\Sigma-2c_o\int_\Omega \frac{1}{z^2}\,dV\,,
\end{equation}
where $\Omega\subset{\bf R}^3$ is the algebraic volume enclosed by the surface $\Sigma$.

The part of $\Sigma$ contained in ${\bf R}^3_+$ can be regarded, after a conformal change of metric, as a surface in the hyperbolic space\footnote{Throughout this paper we will consider the upper half-space model for ${\bf H}^3$ although some of our illustrations will be shown in the ball model.} ${\bf H}^3$. In this setting, the functional \eqref{G} can be expressed as
\begin{equation}\label{GH}
	\mathcal{G}[\Sigma]=\widehat{\mathcal{A}}[\Sigma]-2c_o\widehat{\mathcal{U}}[\Sigma]\,,
\end{equation}
where $\widehat{\mathcal{A}}$ denotes the hyperbolic area and $\widehat{\mathcal{U}}$ is the hyperbolic gravitational potential energy with respect to the constant downward pointing gravitational field $-E_3$, which is given by
\begin{equation}\label{U}
	\widehat{\mathcal{U}}[\Sigma]:=\int_{\widehat{\Omega}} z\,d\widehat{V}\,.
\end{equation}
Here, $\widehat{\Omega}\subset{\bf H}^3$ denotes the (hyperbolic) volume enclosed by $\Sigma$ and $d\widehat{V}$ is the volume element on ${\bf H}^3$. 

Notice that the mean curvature of $\Sigma$ as a surface in ${\bf H}^3$ is given in terms of the objects in the Euclidean space by
\begin{equation}\label{HH}
	\widehat{H}=H\,z+\nu_3\,.
\end{equation}
Hence, the reduced membrane equation \eqref{RME} becomes
\begin{equation}\label{RMEH}
	\widehat{H}+c_o\,z=0\,.
\end{equation}

Surfaces in the Euclidean space ${\bf R}^3$ whose mean curvature depends linearly on their height are used to model liquid drops in a vertical gravitational field \cite{F,W}. When the surfaces are axially symmetric and have their boundaries in a horizontal plane, they are called pendant drops if they lie below this plane, while if they lie above they are referred to as sessile drops.

\section{Regularized Functional $\mathcal{G}_R$}

From now on, we will assume that the interior of $\Sigma$ is embedded in ${\bf R}^3_+$, while its boundary $\partial\Sigma$ lies in the plane $\{z=0\}$. We will further assume that $\Sigma$ meets the plane $\{z=0\}$ orthogonally. Such a surface can be understood, after a conformal change of metric, as a surface in ${\bf H}^3$ approaching the ideal boundary. In this case, both the hyperbolic area $\widehat{\mathcal{A}}$ and the hyperbolic gravitational potential energy $\widehat{\mathcal{U}}$ of $\Sigma$ become infinite and so we will consider their corresponding regularizations.

Let $\underline{z}>0$ and denote, for convenience, $\underline{\Sigma}=\Sigma\cap\{z\geq \underline{z}\}$. We then define the functional $\mathcal{A}_R$ as
\begin{equation}\label{AR}
	\mathcal{A}_R[\Sigma]:=\lim_{\underline{z}\to0^+}\left(\int_{\underline{\Sigma}}\frac{1}{z^2}\,d\Sigma+\frac{1}{\underline{z}}\oint_{\partial\underline{\Sigma}} \star dz\right),
\end{equation}
where $\star dz$ represents the codifferential of $z$. In other words, $\star dz=\partial_n z$, where $n$ is the outward pointing conormal to $\partial\underline{\Sigma}$ and $\partial_n$ represents the derivative in this conormal direction.

\begin{remark} For each $\underline{z}>0$, the area of $\underline{\Sigma}$ is finite and has an expansion in terms of $\underline{z}$. The constant term in that expansion is the Hadamard regularization of the area, known as the renormalized area (see, for instance, \cite{AM,GW}). The functional $\mathcal{A}_R$ as given in \eqref{AR} coincides with the renormalized area (cf., (2.7) of \cite{AM}).
\end{remark}

We also introduce the functional 
\begin{equation}\label{UR}
	\mathcal{U}_R[\Sigma]:=\lim_{\underline{z}\to 0^+}\left(\int_{\underline{\Omega}}\frac{1}{z^2}\,dV-\frac{\lvert \Omega\cap\{z=\underline{z}\}\rvert}{\underline{z}}\right),
\end{equation}
where $\Omega\subset{\bf R}^3$ is the region enclosed by the surface $\Sigma$ and the plane $\{z=0\}$, $\underline{\Omega}=\Omega\cap\{z\geq \underline{z}\}$, and $\lvert \Omega\cap\{z=\underline{z}\}\rvert$ represents the Euclidean area of $\Omega\cap\{z=\underline{z}\}$. The functional $\mathcal{U}_R$ is referred to as the regularization of $\widehat{\mathcal{U}}$. This functional can be rewritten as
\begin{equation}\label{URNew}
	\mathcal{U}_R[\Sigma]=-\int_\Sigma \frac{\nu_3}{z}\,d\Sigma\,.
\end{equation}
Indeed, from the Divergence Theorem we obtain
\begin{eqnarray*}
	\int_{\underline{\Omega}} \frac{1}{z^2}\,dV&=&\int_{\underline{\Omega}} -\nabla'\cdot\left(\frac{E_3}{z}\right)dV =-\int_{\underline{\Sigma}}\frac{\nu_3}{z}\,d\Sigma-\int_{\Omega\cap\{z=\underline{z}\}}\frac{-1}{z}\,d\Sigma\\
	&=&-\int_{\underline{\Sigma}}\frac{\nu_3}{z}\,d\Sigma+\frac{\lvert \Omega\cap\{z=\underline{z}\}\rvert}{\underline{z}}\,.
\end{eqnarray*}
Here, we are denoting by $\nabla'\cdot $ the divergence operator on ${\bf R}^3$. Rearranging and taking the limit as $\underline{z}\to 0^+$ we deduce the expression \eqref{URNew} for $\mathcal{U}_R$.

Therefore, we consider the regularized functional $\mathcal{G}_R$ given by
\begin{equation}\label{GR}
	\mathcal{G}_R[\Sigma]:=\mathcal{A}_R[\Sigma]-2c_o\mathcal{U}_R[\Sigma]\,.
\end{equation}
In the following result we obtain a relation between the regularized functional $\mathcal{G}_R$ and the Helfrich functional $\mathcal{H}_{1,c_o,0}$.

\begin{theorem}\label{thm1} Let $\Sigma$ be a compact surface whose boundary lies in the plane $\{z=0\}$ and assume that $\Sigma$ meets this plane orthogonally. The regularized functional $\mathcal{G}_R$ of $\Sigma$ satisfies
	\begin{equation}\label{equality}
		-\mathcal{G}_R[\Sigma]+\int_\Sigma \left({\widehat H}+c_o\,z\right)^2d\widehat{\Sigma}=\mathcal{H}_{1,c_o,0}[\Sigma]\,,
	\end{equation}
where $\widehat{H}$ is the mean curvature of $\Sigma$ regarded as a surface in ${\bf H}^3$, $d\widehat{\Sigma}$ is the area element on $\Sigma$ induced by the metric on ${\bf H}^3$, and $\mathcal{H}_{1,c_o,0}$ is the Helfrich energy \eqref{H}. In particular, it follows that
\begin{equation}\label{inequality}
	-\mathcal{G}_R[\Sigma]\leq\mathcal{H}_{1,c_o,0}[\Sigma]\,,
\end{equation}
with equality if and only if the reduced membrane equation \eqref{RME} holds on $\Sigma$.
\end{theorem}
{\it Proof.\:} Let $\underline{z}>0$ be sufficiently small and consider the part of $\Sigma$, $\underline{\Sigma}=\Sigma\cap\{z\geq\underline{z}\}$. On $\underline{\Sigma}$ we compute
\begin{eqnarray*}
	\Delta\left(\ln z\right)&=&\nabla\cdot\nabla\left(\ln z\right)=\nabla\cdot\left(\frac{\nabla z}{z}\right)=\frac{\Delta z}{z}-\frac{\lVert \nabla z\rVert^2}{z^2}=2H\frac{\nu_3}{z}-\frac{1-\nu_3^2}{z^2}\\
	&=&2H\frac{\nu_3}{z}+\frac{\nu_3^2}{z^2}-\frac{1}{z^2}\,,
\end{eqnarray*}
where we have used that $\Delta z=2H\nu_3$. (Here, $\nabla\cdot$ denotes the divergence operator on the surface $\underline{\Sigma}$.)

Integrating over $\underline{\Sigma}$ and applying the Divergence Theorem to the left-hand side, we get
$$\frac{1}{\underline{z}}\oint_{\partial\underline{\Sigma}} \star dz=\int_{\underline{\Sigma}}\left(2H\frac{\nu_3}{z}+\frac{\nu_3^2}{z^2}\right)d\Sigma-\int_{\underline{\Sigma}}\frac{1}{z^2}\,d\Sigma\,.$$
Combining this with the expression of $\mathcal{U}_R$ given in \eqref{URNew}, we obtain
\begin{eqnarray*}
	\int_{\underline{\Sigma}}\frac{1}{z^2}\,d\Sigma+\frac{1}{\underline{z}}\oint_{\partial\underline{\Sigma}}\star dz-2c_o \mathcal{U}_R\left[\,\underline{\Sigma}\,\right]=\int_{\underline{\Sigma}}\left(2\left[H+c_o\right]\frac{\nu_3}{z}+\frac{\nu_3^2}{z^2}\right)d\Sigma\\
	=\int_{\underline{\Sigma}}\left(H+c_o+\frac{\nu_3}{z}\right)^2d\Sigma-\int_{\underline{\Sigma}}\left(H+c_o\right)^2d\Sigma\,.
\end{eqnarray*}
By letting $\underline{z}\to 0^+$, we conclude with
\begin{equation*}
	\mathcal{G}_R[\Sigma]=\mathcal{A}_R[\Sigma]-2c_o\mathcal{U}_R[\Sigma]=\int_\Sigma\left(H+c_o+\frac{\nu_3}{z}\right)^2d\Sigma-\mathcal{H}_{1,c_o,0}[\Sigma]\,,
\end{equation*}
where we have used the definition of the Helfrich energy given in \eqref{H}. Substituting \eqref{HH} in the integral above, the first statement follows.

Finally, since the first term on the right-hand side of the previous equation is non negative, the inequality \eqref{inequality} holds with equality if and only if \eqref{RME} holds. {\bf q.e.d.}
\\

We finish this section by illustrating the equality \eqref{equality} for the special case of the round hemisphere. Let $c_o>0$ and denote by $S^2_+$ the hemisphere of radius $R=1/c_o$ in ${\bf R}^3_+$ having its boundary in the plane $\{z=0\}$. Note that $H+c_o\equiv 0$ holds. Using the spherical coordinates $X(\theta,\phi)=(R\cos\theta\sin\phi,R\sin\theta\sin\phi,R\cos\phi)$ to parameterize $S_+^2$, we compute
$$\int_{\underline{S}_+^2}\frac{1}{z^2}\,d\Sigma+\frac{1}{\underline{z}}\oint_{\partial\underline{S}_+^2}\star dz=2\pi\int_0^{\underline{\phi}}\frac{\sin\phi}{\cos^2\phi}\,d\phi-2\pi\frac{\sin\underline{\phi}}{\cos\underline{\phi}}=2\pi\left(\frac{1-\sin\underline{\phi}}{\cos\underline{\phi}}-1\right),$$
where $\underline{\phi}\in(0,\pi/2)$ is such that $\underline{z}=R\cos\underline{\phi}$. Observe that as $\underline{z}\to 0^+$, $\underline{\phi}\to\pi/2^-$. Hence, taking this limit, we deduce from \eqref{AR} that
\begin{equation}\label{ARS}
	\mathcal{A}_R[S_+^2]=\lim_{\underline{\phi}\to \pi/2^-} 2\pi\left(\frac{1-\sin\underline{\phi}}{\cos\underline{\phi}}-1\right)=-2\pi\,.
\end{equation}
Using \eqref{URNew} we also compute
\begin{equation}\label{URS}
	\mathcal{U}_R[S_+^2]=-2\pi\int_0^{\pi/2} R\sin\phi\,d\phi=-2\pi R\,,
\end{equation}
since $\nu_3=\cos\phi=z/R$ holds. Combining \eqref{ARS} and \eqref{URS}, we deduce that
\begin{equation}\label{GRS}
	\mathcal{G}_R[S_+^2]=\mathcal{A}_R[S_+^2]-2c_o\mathcal{U}_R[S_+^2]=-2\pi+4\pi=2\pi\,,
\end{equation}
because $c_oR=1$. On the other hand, one can compute $\mathcal{G}_R[S_+^2]$ using the right-hand side of \eqref{equality}. Since $H+c_o\equiv 0$ holds on $S_+^2$, it follows that $\mathcal{H}_{1,c_o,0}[S_+^2]=0$. In addition, $S_+^2$ is minimal in ${\bf H}^3$ and so $\widehat{H}\equiv 0$ holds. Then, \eqref{equality} simplifies to
$$\mathcal{G}_R[S_+^2]=\int_{S_+^2} c_o^2z^2\,d\widehat{\Sigma}=\int_{S_+^2} c_o^2\,d\Sigma=2c_o^2\int_0^{\pi/2} R^2\sin\phi\,d\phi=2\pi\,,$$
as expected.

\section{Critica of $-\mathcal{G}_R$}

Throughout this section, we will define admissible surfaces to be those surfaces $\Sigma$ contained in the upper half-space ${\bf R}^3_+$ which meet the plane $\{z=0\}$ at their boundaries $\partial\Sigma$ in such a way that $\nu_3/z$ smoothly extends to $z=0$. In particular, $\nu_3=0$ must hold along the boundary of an admissible surface.

An admissible variation will be a variation through admissible surfaces. This means that, for admissible variations, the associated variational field must be tangent to the plane $\{z=0\}$ and that the orthogonality of the intersection of $\Sigma$ with $\{z=0\}$ must be preserved through the variation. Let $\psi\in\mathcal{C}^\infty(\Sigma)$ be a smooth function defined on $\Sigma$ and consider an arbitrary variational vector field $\delta X=\psi\nu+(\delta X)^T$. To have an admissible variation, we must impose the following restrictions on $\delta X$: on the one hand, the vector field $\delta X$ is tangent to $\{z=0\}$ if and only if $\delta X\cdot E_3=0$ along $\partial\Sigma$. Since along $\partial\Sigma$, the outward pointing conormal $n$ satisfies $n=-E_3$, the tangential condition translates to
$$\delta X\cdot E_3=\psi\nu_3-(\delta X)^T\cdot n=-(\delta X)^T\cdot n=0\,,$$
along the boundary. On the other hand, the condition that the variation preserves $\nu_3=0$ along $\partial\Sigma$ gives
$$\delta\nu_3=-\nabla\psi\cdot E_3+d\nu(\delta X^T)\cdot E_3=\left(\nabla\psi-d\nu(\delta X^T)\right)\cdot n=0\,,$$
along the boundary. In conclusion, any normal variation $\delta X=\psi\nu$ is admissible if and only if $\partial_n\psi=0$ on $\partial\Sigma$.

We will next obtain the equations characterizing equilibria for the functional\footnote{From Theorem \ref{thm1}, minimization of the functional $-\mathcal{G}_R$ is closely related to the minimization of the Helfrich energy. Hence, we will consider the functional $-\mathcal{G}_R$ rather than the one with the positive sign. Nonetheless, equilibria are the same regardless of the choice of sign.} $-\mathcal{G}_R$ defined in \eqref{GR}.

\begin{theorem}\label{ELequations} An admissible surface $\Sigma$ is in equilibrium for the functional $-\mathcal{G}_R$ if and only if the reduced membrane equation \eqref{RME} holds on the interior of $\Sigma$ and $\partial_nH=0$ holds along the boundary $\partial\Sigma$. (Recall that $\partial_n$ denotes the derivative in the conormal direction.)
	
Further, for any equilibrium surface
$$2c_o\mathcal{U}_R[\Sigma]=2c_o\int_\Sigma\left(H+c_o\right)d\Sigma=0\,,$$
holds.
\end{theorem}
{\it Proof.\:} Compactly supported normal variations can be applied to deduce that the reduced membrane equation \eqref{RME} holds on the interior for any equilibrium. We refer the reader to Theorem 4.2 of \cite{PP2} for details.

Assume now that the reduced membrane equation \eqref{RME} holds. To obtain the boundary condition we will next compute the first variation formula of $-\mathcal{G}_R$ employing the equality \eqref{equality}. Since \eqref{RME} holds (cf., \eqref{RMEH}), the integral term in \eqref{equality} vanishes and, for arbitrary variations $\delta X=\psi\nu+(\delta X)^T$, we have
\begin{eqnarray*}
	\delta\left(-\mathcal{G}_R\right)[\Sigma]&=&\delta\mathcal{H}_{1,c_o,0}[\Sigma]=\int_\Sigma \mathcal{L}[H+c_o]\psi\,d\Sigma\\&&+\oint_{\partial\Sigma}\left(\left[H+c_o\right]\partial_n\psi-\partial_nH\psi+\left[H+c_o\right]^2(\delta X)^T\cdot n\right)ds\,,
\end{eqnarray*}
where $\mathcal{L}$ is the second order elliptic operator defined by
$$\mathcal{L}[f]:=\Delta f+2\left(H[H-c_o]-K\right)f\,.$$
Details about the computation of the first variation formula for the Helfrich energy $\mathcal{H}_{1,c_o,0}$ can be found in Appendix A of \cite{PP1}.

From Proposition 4.1 of \cite{PP2}, surfaces satisfying the reduced membrane equation \eqref{RME} also satisfy \eqref{ELH} and, hence, $\mathcal{L}[H+c_o]=0$ holds on $\Sigma$. Therefore, the surface integral above vanishes. Moreover, recall that for admissible variations $(\delta X)^T\cdot n=0=\partial_n\psi$ hold along $\partial\Sigma$. We then conclude with
$$\delta\left(-\mathcal{G}_R\right)[\Sigma]=-\oint_{\partial\Sigma}\partial_nH\,\psi\,ds\,.$$
Hence, if $\delta(-\mathcal{G}_R)[\Sigma]=0$ holds for every admissible variation, $\partial_nH=0$ must hold on $\partial\Sigma$. 

Finally, we will prove that the regularization $\mathcal{U}_R$ of the (hyperbolic) gravitational potential energy $\widehat{\mathcal{U}}$ at an equilibrium surface $\Sigma$ satisfies $2c_o\mathcal{U}_R[\Sigma]=0$. Since for an equilibrium surface the reduce membrane equation \eqref{RME} holds, it follows from \eqref{URNew} that
\begin{equation}\label{URHelp}
	\mathcal{U}_R[\Sigma]=-\int_\Sigma \frac{\nu_3}{z}\,d\Sigma=\int_\Sigma\left(H+c_o\right)d\Sigma\,.
\end{equation} 
In \cite{PP2}, using variations $\delta X=X$, it was shown that for any Helfrich surface (which is the case of $\Sigma$ since \eqref{RME} holds),
$$2c_o\int_\Sigma \left(H+c_o\right)d\Sigma=\oint_{\partial\Sigma}\left([H+c_o]\partial_nq-q\partial_nH+[H+c_o]^2\partial_n X^2/2\right)ds\,,$$
holds, where $q:=X\cdot\nu$ is the support function. We claim that the boundary integral on the right-hand side vanishes. To see this, observe that since $\nu_3\equiv 0$ holds along $\partial\Sigma$, the surface $\Sigma$ meets the plane $\{z=0\}$ in a right angle along $\partial\Sigma$. Therefore, the outward pointing unit conormal to $\partial\Sigma$ satisfies $n=-E_3$. We then have (see (2) of \cite{PP1})
$$0=-\partial_sE_3=n_s=-\kappa_gT+\tau_g\nu\,,$$
where $s$ denotes the arc length parameter along $\partial\Sigma$, $T$ is the unit tangent vector to $\partial\Sigma$, $\kappa_g$ is the geodesic curvature of $\partial\Sigma$, and $\tau_g$ is the geodesic torsion of $\partial\Sigma$. Hence, we obtain that $\kappa_g\equiv 0$ and $\tau_g\equiv 0$, which mean that each boundary component of $\partial\Sigma$ is a geodesic and line of curvature, respectively. By the latter statement, we can write $\nu_s=-\kappa_n T$ (see again (2) of \cite{PP1}) and, hence, the normal curvature $\kappa_n$ is a principal curvature, say $\kappa_2$ along the boundary. Therefore,
$$\partial_nq=\partial_n\left(X\cdot\nu\right)=n\cdot\nu+X\cdot\partial_n\nu=X\cdot d\nu(n)=-(2H-\kappa_n)X\cdot n=0\,,$$
where the last equality follows since $\tau_g\equiv 0$ holds along $\partial\Sigma$. (Recall that $\kappa_1=2H-\kappa_n$ along $\partial\Sigma$). In addition, we have the boundary condition $\partial_nH=0$ and
$$\partial_n(X^2/2)=X\cdot \partial_n X=X\cdot n=0\,.$$
Combining these three things we deduce that the boundary integral above vanishes. Thus, together with \eqref{URHelp} we finish the proof. {\bf q.e.d.}
\\

The reduced membrane equation \eqref{RME} in combination with the boundary condition $\partial_nH=0$ and the smooth extension of an equilibrium surface $\Sigma$ for $-\mathcal{G}_R$ to $z=0$ describe the boundary $\partial\Sigma$ of such surfaces.

\begin{prop}\label{propboundary} Let $\Sigma$ be an equilibrium surface for the functional $-\mathcal{G}_R$. The boundary components of $\Sigma$ are geodesic lines of curvature along which
\begin{equation}\label{k1k2}
	\kappa_1-\kappa_2=2c_o\,,
\end{equation}
holds, where $\kappa_1$ and $\kappa_2$ are the principal curvatures of $\Sigma$.

Moreover, if $c_o\neq 0$,
\begin{equation*}
	\lim_{\underline{z}\to 0^+} \frac{\kappa_1-\kappa_2-2c_o}{\underline{z}}=0\,.
\end{equation*}
\end{prop}
{\it Proof.\:} Consider a surface $\Sigma$ satisfying the reduced membrane equation \eqref{RME} and which extends smoothly to $z=0$. As shown in the last part of the proof of Theorem \ref{ELequations}, the boundary components of $\Sigma$ are geodesic lines of curvature. Hence, it only remains to prove the relations between the principal curvatures at the boundary.

A straightforward calculation using \eqref{RME} shows that for any point in $\Sigma$ with $z>0$, there holds
$$\nabla H=\nabla\left(-\frac{\nu_3}{z}\right)=-\frac{\nabla\nu_3}{z}+\frac{\nu_3\nabla z}{z^2}=\frac{-\nabla\nu_3-(H+c_o)\nabla z}{z}\,.$$
In the last equality, we have used once again that the reduced membrane equation \eqref{RME} holds on $\Sigma$. Observe that along any boundary component
$$\nabla\nu_3=\nabla\nu\cdot E_3=d\nu(E_3)=d\nu(\nabla z)=-\kappa_1\nabla z\,,$$
since the boundary component is a line of curvature. Recall that $\kappa_1$ is the principal curvature such that $\kappa_1=2H-\kappa_n$ along the boundary.

From the above computation and since $\Sigma$ extends smoothly to $z=0$, we deduce that at $z=0$ the numerator of $\nabla H$ must vanish. Thus,
$$0=-\nabla\nu_3-(H+c_o)\nabla z=\left(\kappa_1-(H+c_o)\right)\nabla z=-(\kappa_1-H-c_o)E_3\,,$$
where the last equality holds since $\nabla z=-E_3$ along the boundary. Therefore, $\kappa_1-H=c_o$ along $\partial\Sigma$, or equivalently $\kappa_1-\kappa_2=2c_o$.

For the second identity, since $c_o\neq 0$, it follows from $\kappa_1-\kappa_2=2c_o$ that there are no umbilic points at $z=0$. Hence, near any boundary component we can pick up an orthonormal frame of principal directions $\{e_1,e_2\}$, such that $e_1=n$ and $e_2=T$ along $z=0$ (the latter is possible since boundary components are lines of curvature). In this frame, at $z=\underline{z}>0$ we get from the above computation that
\begin{eqnarray*}
	e_1(H)&=&\nabla H\cdot e_1=\frac{-\nabla\nu_3\cdot e_1-(H+c_o)\nabla z\cdot e_1}{\underline{z}}=\frac{\kappa_1-H+c_o}{\underline{z}}\,e_1(z)\\
	&=&\frac{\kappa_1-\kappa_2-2c_o}{2\underline{z}}\,e_1(z)\,,
\end{eqnarray*}
since $\nabla\nu_3\cdot e_1=e_1(\nu_3)=d\nu(e_1)\cdot E_3=d\nu(e_1)\cdot \nabla z=-\kappa_1 e_1\cdot\nabla z$ and $\nabla z\cdot e_1=e_1(z)$. Taking the limit when $\underline{z}\to 0^+$ and using the boundary condition $\partial_nH=0$, we get
$$0=\partial_nH=e_1(H)=\frac{-1}{2}\lim_{\underline{z}\to 0^+} \frac{\kappa_1-\kappa_2-2c_o}{\underline{z}}\,,$$
since $e_1(z)=\partial_n z=-1$ at $z=0$. {\bf q.e.d.}
\\

In the following theorems we will show a couple of rigidity results for equilibrium surfaces for $-\mathcal{G}_R$.

\begin{theorem}\label{newc=0} Genus zero equilibria for the functional $-\mathcal{G}_R$ with $c_o=0$ are hemispheres $S_+^2$ of arbitrary radii having their boundaries in the plane $\{z=0\}$ (ie., totally geodesic discs of ${\bf H}^3$).
\end{theorem}
{\it Proof.\:} Let $\Sigma$ be an equilibrium surface for $-\mathcal{G}_R$ with $c_o=0$ and assume that $\Sigma$ has genus zero and an arbitrary number $\mathfrak{b}$ of boundary components. On $\Sigma$ we consider the quadratic differential $\Psi=\Phi/z$, where $\Phi$ is the classical (Euclidean) quadratic Hopf differential. We will prove that, since $\partial_nH=0$ along $\partial\Sigma$, $\Psi$ extends holomorphically to a closed doubled surface of genus $\mathfrak{g}=\mathfrak{b}-1$. The proof of the statement then will follow by showing that $\Psi$ has enough zeros to force it to be identically zero.

Let $C$ be any connected component of $\partial\Sigma$ and consider, near $C$, conformal coordinates $w=u+iv$ adapted to $C$. That is, take $v=0$ to be $C$ and $v\geq 0$ the interior of the surface. Hence, $u$ is the direction of $C$ and $v$ is the direction orthogonal to $C$ pointing inward. The metric of $\Sigma$ is then given by $E(du^2+dv^2)$. We denote by $Ldu^2+2Mdudv+Ndv^2$ the second fundamental form of $\Sigma$ in these local coordinates.
	
Let $\Phi$ be the (Euclidean) quadratic Hopf differential \cite{Hopf} and construct the quadratic differential
$$\Psi:=\frac{\Phi}{z}\,.$$
\emph{After a conformal change of metric, the quadratic differential $\Psi$ is the hyperbolic Hopf differential, which for minimal surfaces in ${\bf H}^3$ (see \eqref{RMEH} for $c_o=0$) is holomorphic.} In the local coordinates $w=u+iv$, the differential $\Psi$ has coefficient
$$q:=\frac{Q}{z}=\frac{L-N}{2z}-i\frac{M}{z}\,,$$
where $Q$ is the coefficient of the Hopf differential $\Phi$. Since $c_o=0$, it follows from \eqref{k1k2} that the boundary component $C$ is composed by umbilic points. These points are the zeros of the Hopf differential so, along $C$, $Q\equiv 0$ holds. Hence, $q$ is not singular at $C$.

We next compute the value of $q$ at the boundary component $C$. For that, we employ L'Hopital's rule to obtain
$$q(u,0)=\frac{Q_v(u,0)}{z_v(u,0)}=\frac{L_v-N_v}{2z_v}-i\frac{M_v}{z_v}\,.$$
Observe that
\begin{equation}\label{zv}
	z_v(u,0)=X_v(u,0)\cdot E_3=\sqrt{E}\,,
\end{equation}
since $X_v(u,0)$ is parallel to $E_3$ along $C$, points in the same direction, and has magnitude $\sqrt{E}$. To further simplify $q(u,0)$ we use the Codazzi equations of $\Sigma$. These equations can be described in terms of the derivative of $Q$ as $Q_{\overline{w}}=EH_w$. Taking the real and imaginary parts of this equation, they read
\begin{eqnarray}
	\frac{L_u-N_u}{2}+M_v&=&EH_u\,,\label{C1}\\
	\frac{L_v-N_v}{2}-M_u&=&-EH_v\,.\label{C2}
\end{eqnarray}
Recall that $C$ is a line of curvature so $M(u,0)=0$. Differentiating this along $C$, we get $M_u(u,0)=0$ as well. In addition, it follows that at $v=0$, the coordinate directions are principal directions. Hence, 
$$\kappa_1(u,0)=\frac{N(u,0)}{E(u,0)}\,,\quad\quad\quad\kappa_2(u,0)=\frac{L(u,0)}{E(u,0)}\,,$$
and \eqref{k1k2} reduces to $N(u,0)-L(u,0)=0$. The derivative along $C$ of this identity gives $N_u(u,0)=L_u(u,0)$. Using these derivatives the Codazzi equations \eqref{C1}-\eqref{C2} simplify to
\begin{eqnarray*}
	M_v&=&EH_u\,,\\
	L_v-N_v&=&-2EH_v\,,
\end{eqnarray*}
along $C$. We then conclude from this and \eqref{zv} that
\begin{equation}\label{q(u,0)}
	q(u,0)=-\sqrt{E}H_v-i\sqrt{E}H_u=-i\sqrt{E}H_u\,,
\end{equation}
where the last equality follows since for an equilibrium surface for $-\mathcal{G}_R$, $0=\partial_nH=-H_v$ holds along the boundary.

We now verify that $\Psi$ is holomorphic on $\Sigma$. For this purpose, we compute
\begin{equation}\label{bar}
	q_{\overline{w}}=\frac{Q_{\overline{w}}}{z}-\frac{Qz_{\overline{w}}}{z^2}=\frac{EH_w}{z}-\frac{Qz_{\overline{w}}}{z^2}\,,
\end{equation}
where in the last equality we have employed the Codazzi equations in the complex setting. Differentiating the reduced membrane equation \eqref{RME} for $c_o=0$, we obtain
$$H_w=-\left(\frac{\nu_3}{z}\right)_w=-\frac{(\nu_3)_w}{z}+\frac{\nu_3 z_w}{z^2}=-\frac{(\nu_3)_w}{z}-H\frac{z_w}{z}=\frac{Qz_{\overline{w}}}{Ez}\,,$$
since $(\nu_3)_w=-Hz_w-Qz_{\overline{w}}/E$ holds. Therefore, substituting this in \eqref{bar}, we deduce that $q_{\overline{w}}=0$ and so $\Psi$ is holomorphic. Since the real part of $q$ at $C$ vanishes (see \eqref{q(u,0)}) and $q$ is not singular at $C$, the holomorphic quadratic differential $\Psi$ can be holomorphically reflected to the closed doubled genus $\mathfrak{g}=\mathfrak{b}-1$ surface.
 
A nonzero holomorphic quadratic differential on a closed surface of genus $\mathfrak{g}=\mathfrak{b}-1$ has $4\mathfrak{g}-4=4\mathfrak{b}-8$ zeros, counted with multiplicity. Observe that, from Proposition \ref{propboundary}, since each boundary component is a geodesic of $\Sigma$, its curvature $\kappa$ (as a curve in the plane $\{z=0\}$) satisfies
$$\kappa^2=\kappa_g^2+\kappa_n^2=\kappa_n^2\,,$$
where the normal curvature $\kappa_n=\kappa_2$ is a principal curvature along the boundary component. By the Four-Vertex Theorem, the curvature $\kappa$ of each boundary component has, at least, four vertices and so at each boundary component there are, at least, four points at which $(\kappa_n)_u=(\kappa_2)_u=0$. Furthermore, from \eqref{k1k2}, we get $H=\kappa_2$ along each boundary component. Thus, the vertices of the boundary component correspond with points where $H_u=0$, that is, zeros of the extended $\Psi$ (see \eqref{q(u,0)}). Therefore, we have shown the existence of, at least, $4\mathfrak{b}$ zeros of the holomorphic quadratic differential extension of $\Psi$, reaching a contradiction unless this extension is zero and so $\Psi\equiv 0$ identically.

Finally, from $\Psi\equiv 0$, we deduce that on the original surface $\Sigma$, the quadratic Hopf differential $\Phi\equiv 0$ identically and so the surface is totally umbilical. This finishes the proof. {\bf q.e.d.}

\begin{remark} In the specific case $c_o=0$, the functional $-\mathcal{G}_R\equiv-\mathcal{A}_R$ is just the renormalized area functional. In Theorem 7.1 of \cite{AM}, employing a refined version of the asymptotic maximum principle, it was shown that an equilibrium for $-\mathcal{A}_R$ with connected boundary is a totally geodesic disc of ${\bf H}^3$ (see also Remark 7.2 of \cite{AM}). The proof of Theorem \ref{newc=0} above is essentially different to that of \cite{AM} and it applies to a different class of surfaces, namely, genus zero surfaces with an arbitrary number of boundary components.
\end{remark}

We next show that equilibrium surfaces for $-\mathcal{G}_R$ with a circular boundary component are axially symmetric. In \cite{W}, Wente showed the axial symmetry of {\it all embedded} pendant and sessile drops in ${\bf R}^3$ using Alexandrov's reflection method. Our approach here will use that surfaces satisfying the reduced membrane equation \eqref{RME} also satisfy \eqref{ELH} and the ideas of Lemma 4.1 of \cite{PP2}.

\begin{theorem}\label{axsym}
	Equilibria for the functional $-\mathcal{G}_R$ with a circular boundary component are axially symmetric surfaces whose axis of rotation is the $z$-axis.
\end{theorem}
{\it Proof.\:} Let $\Sigma$ be an equilibrium surface for $-\mathcal{G}_R$. From Theorem \ref{ELequations}, the reduced membrane equation \eqref{RME} holds on $\Sigma$. Hence, as noted above (see Proposition 4.1 of \cite{PP2}), equation \eqref{ELH} is also satisfied on $\Sigma$.

In addition, equilibria for $-\mathcal{G}_R$ also satisfy $\partial_nH=0$ along their boundaries. Moreover, Proposition \ref{propboundary} shows that each boundary component is a geodesic line of curvature along which $H-\kappa_n-c_o=0$ holds. The latter is an immediate consequence of \eqref{k1k2} together with $\kappa_1=2H-\kappa_n$ and $\kappa_2=\kappa_n$ along $\partial\Sigma$.

Let $C$ be the connected component of $\partial\Sigma$ which is a round circle. We will next show that along $C$, $H$ is constant. Since $C$ is a geodesic, its constant curvature $\kappa$ (as a curve in the plane $\{z=0\}$) satisfies $\kappa^2=\kappa_g^2+\kappa_n^2=\kappa_n^2$. Therefore, $\kappa_n$ is constant along $C$ and so is $H=\kappa_n+c_o$.

The remainder of the proof follows by applying Lemma 4.1 of \cite{PP2}. For the sake of clarity, we sketch the ideas here. Let $F$ denote the fourth order self-adjoint operator associated with the second variation of the Helfrich energy $\mathcal{H}_{1,c_o,0}$. That is, for normal variations $\delta X=f\nu$, $f\in\mathcal{C}_o^\infty(\Sigma)$,
$$\delta^2 {\mathcal H}_{1,c_o, 0}[\Sigma]=\int_\Sigma f\,F[f]\:d\Sigma+\frac{1}{2}\oint_{\partial\Sigma} L[f]\partial_nf\,ds\:,$$
where $L$ is the second order operator defined by $L[f]:=\Delta f+\lVert d\nu\rVert^2f$. The explicit computation of the operator $F$ can be found in Section 3.2 of \cite{PP4}, while the above expression of $\delta^2\mathcal{H}_{1,c_o,0}$ is given in Theorem 3.1 of \cite{PP4}.

Since $E_3\times X$ is a Killing vector field on ${\bf R}^3$ which generates a rotation about $E_3$, it is clear that its normal component $\psi:=E_3\times X\cdot \nu$ satisfies $F[\psi]=0$ with $\psi\equiv 0$ along the boundary circle $C$. Calculations given in Lemma 4.1 of \cite{PP2} show that, along $C$, $\partial_n \psi \equiv 0$ as a consequence of $\tau_g\equiv 0$ and $n=-E_3$. In addition, $\partial_n^2\psi\equiv 0$ holds from the previous cases and $H$ being constant along $C$. Finally, $\partial_n^3\psi\equiv0$ holds from the considerations above and the boundary condition $\partial_nH=0$. One then concludes that $\psi\equiv 0$ on all of the surface $\Sigma$ by the Cauchy-Kovalevskaya Theorem. Since $\psi$ is the normal component of the variation obtained by rotating the surface about a vertical axis, this  means that the surface is invariant with respect to this rotation.
{\bf q.e.d.}

\begin{remark} For the case $c_o=0$, Theorem \ref{axsym} shows that equilibria for $-\mathcal{G}_R\equiv-\mathcal{A}_R$ with a circular boundary component are totally geodesic discs of ${\bf H}^3$. Indeed, the other axially symmetric minimal surfaces of ${\bf H}^3$ are (spherical) catenoids, which do not satisfy the boundary condition $\partial_nH=0$.
\end{remark}

Combining the Euler-Lagrange equations associated to the functional $-\mathcal{G}_R$ obtained in Theorem \ref{ELequations} and the result of Theorem \ref{axsym} regarding the axial symmetry of equilibria with a circular boundary component, we conclude with the following description.

\begin{cor} Equilibria for $-\mathcal{G}_R$ with a circular boundary component are axially symmetric surfaces satisfying \eqref{RME} bounded by geodesic lines of curvature parallels located at $\{z=0\}$ along which $\partial_n H=0$ holds.
\end{cor}

We will next briefly explain how to obtain examples of equilibria for $-\mathcal{G}_R$. We will focus our attention on topological discs. Axially symmetric surfaces (whose axis of rotation is the $z$-axis) satisfying the reduced membrane equation \eqref{RME} can be generated by solving the system of first order differential equations
\begin{eqnarray}
	r'(\sigma)&=&\cos\varphi(\sigma)\,,\label{ODE1}\\
	z'(\sigma)&=&\sin\varphi(\sigma)\,,\label{ODE2}\\
	\varphi'(\sigma)&=&-2\frac{\cos\varphi(\sigma)}{z(\sigma)}-\frac{\sin\varphi(\sigma)}{r(\sigma)}-2c_o\,,\label{ODE3}
\end{eqnarray}
where $\gamma(\sigma)=(r(\sigma),z(\sigma))$ is the arc length (measured from the top of the surface) parameterized generating curve and $\varphi(\sigma)$ is the function representing the angle between the positive part of the $r$-axis and the tangent vector to $\gamma$. Here, the derivative with respect to the arc length parameter $\sigma$ is denoted by $\left(\,\right)'$. As noticed in previous papers (see, for instance, \cite{LPP,PP2,PP3}), equation \eqref{ODE3} carries both the positive and negative signs in front of $c_o$. Up to the transformation $z\mapsto -z$, we may assume that this sign is negative. In addition, possibly applying Proposition 3.4 of \cite{LPP}, we will also assume that $c_o\geq 0$ holds.

To obtain topological discs, we need to impose the following initial conditions:
\begin{equation}\label{ic}
	r(0)=0\,,\quad\quad\quad z(0)=z_o\,,\quad\quad\quad\varphi(0)=0\,,
\end{equation}
where $z_o\in{\bf R}\setminus\{0\}$ is considered as a parameter\footnote{Since we are fixing the sign of $c_o$ employing Proposition 3.4 of \cite{LPP}, we need to allow for negative initial heights. The surfaces obtained rotating generating curves for $z_o<0$ will need to be transformed according to $z\mapsto -z$ so as to be contained in ${\bf R}^3_+$.}. 

If $c_o=0$, the unique solution to the system \eqref{ODE1}-\eqref{ODE3} with initial conditions \eqref{ic} is a part of a circle of radius $\lvert z_o\rvert$. For the case $c_o\neq 0$, in Section 3.1 of \cite{LPP} the curves $\gamma(\sigma)=(r(\sigma),z(\sigma))$ obtained from solutions of \eqref{ODE1}-\eqref{ODE3} with initial conditions \eqref{ic} were geometrically described. In particular, it was shown that, when $z_o>-1/c_o$, the curves $\gamma$ meet $\{z=0\}$ orthogonally. In addition to meeting $\{z=0\}$ orthogonally, if one aims to obtain equilibria for $-\mathcal{G}_R$, the boundary condition $\partial_nH=0$ at $z=0$ must also be satisfied. As shown in \cite{LPP}, if $z_o\in(-1/c_o,0)$ this condition cannot hold. Nonetheless, there is numerical evidence that there exists an infinite discrete set of values $z_o>0$ for which $\partial_nH=0$ holds at $z=0$. Recall that $\partial_nH=0$ at $z=0$ is equivalent to $\varphi''=0$ at $z=0$ and compare with Figure 3 of \cite{LPP}.

\begin{figure}[h!]
	\makebox[\textwidth][c]{
		\begin{subfigure}[b]{0.4\linewidth}
			\includegraphics[width=\linewidth]{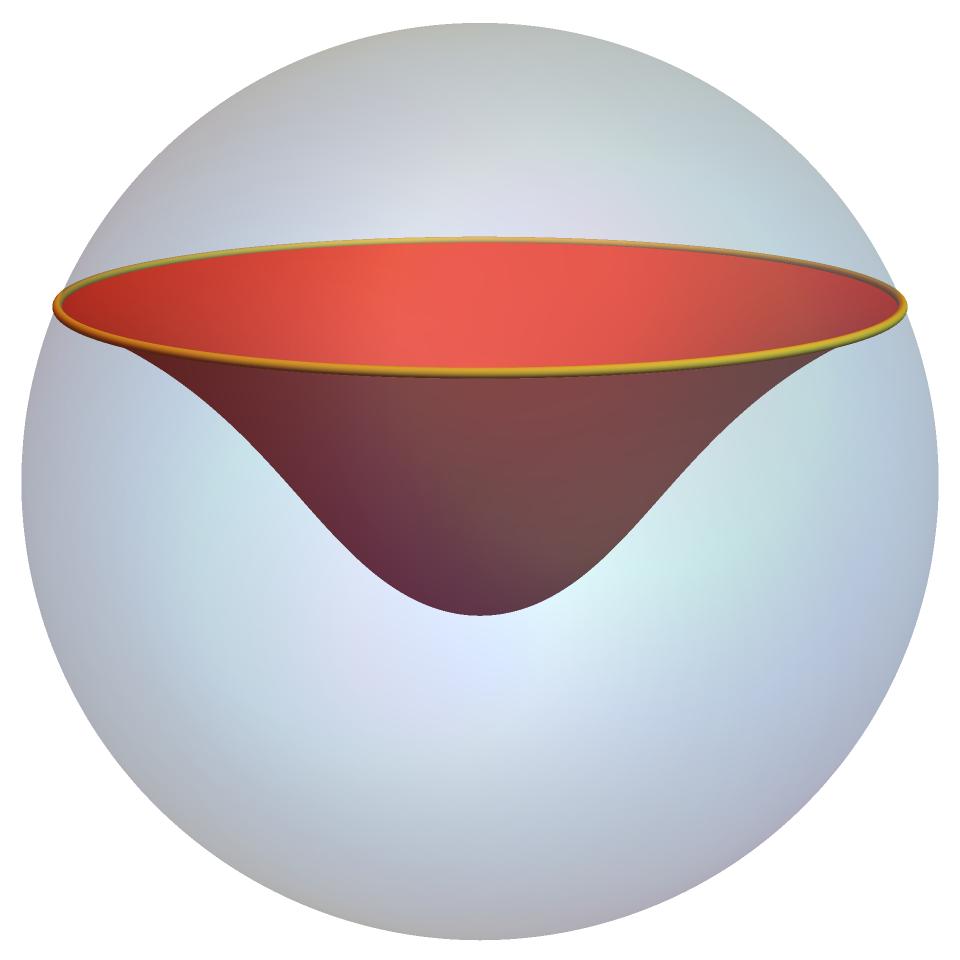}
		\end{subfigure}
		\,\,
		\begin{subfigure}[b]{0.4\linewidth}
			\includegraphics[width=\linewidth]{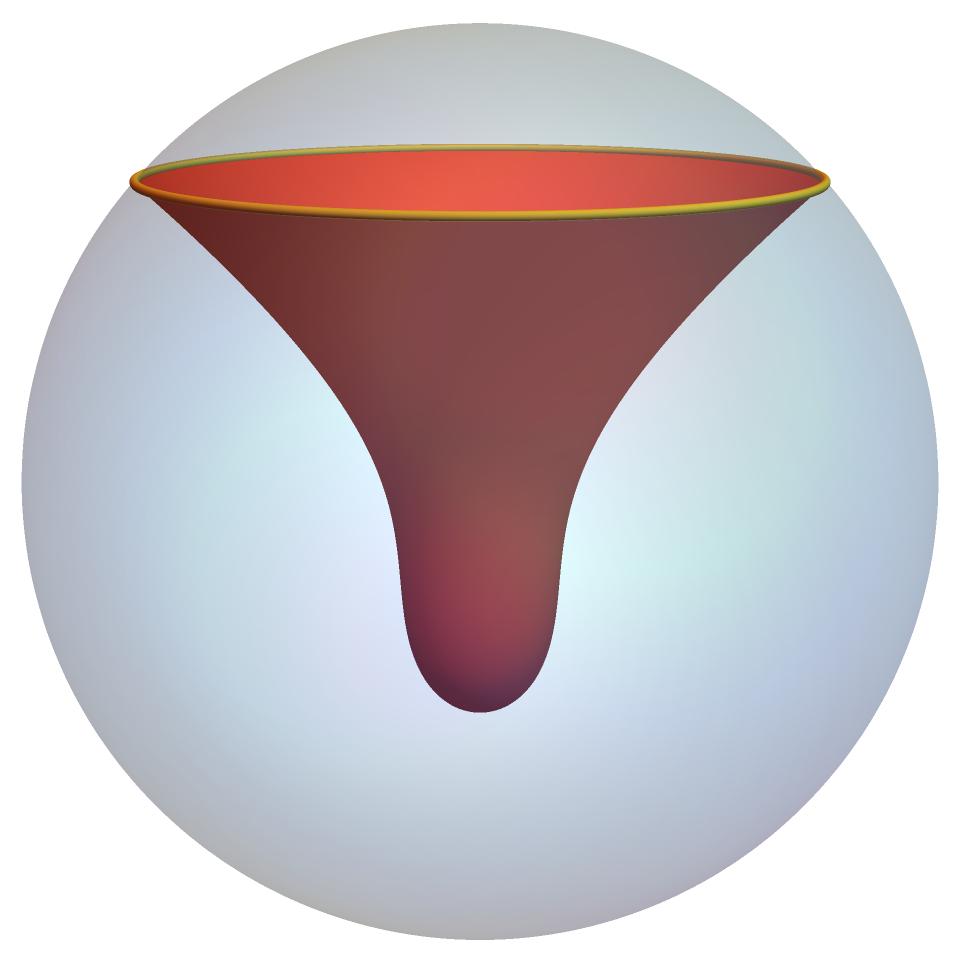}
		\end{subfigure}
		\,\,
		\begin{subfigure}[b]{0.4\linewidth}
			\includegraphics[width=\linewidth]{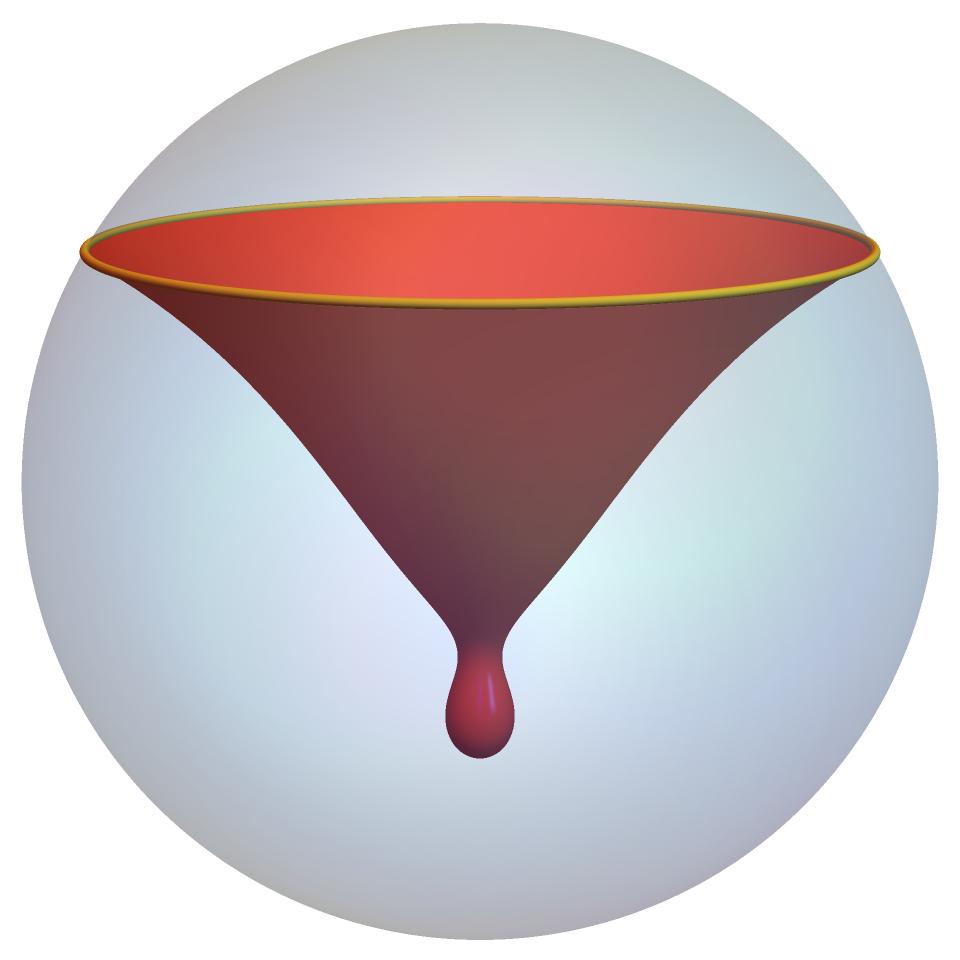}
	\end{subfigure}}\vspace{0.35cm}
	\makebox[\textwidth][c]{
		\begin{subfigure}[b]{0.4\linewidth}
			\includegraphics[width=\linewidth]{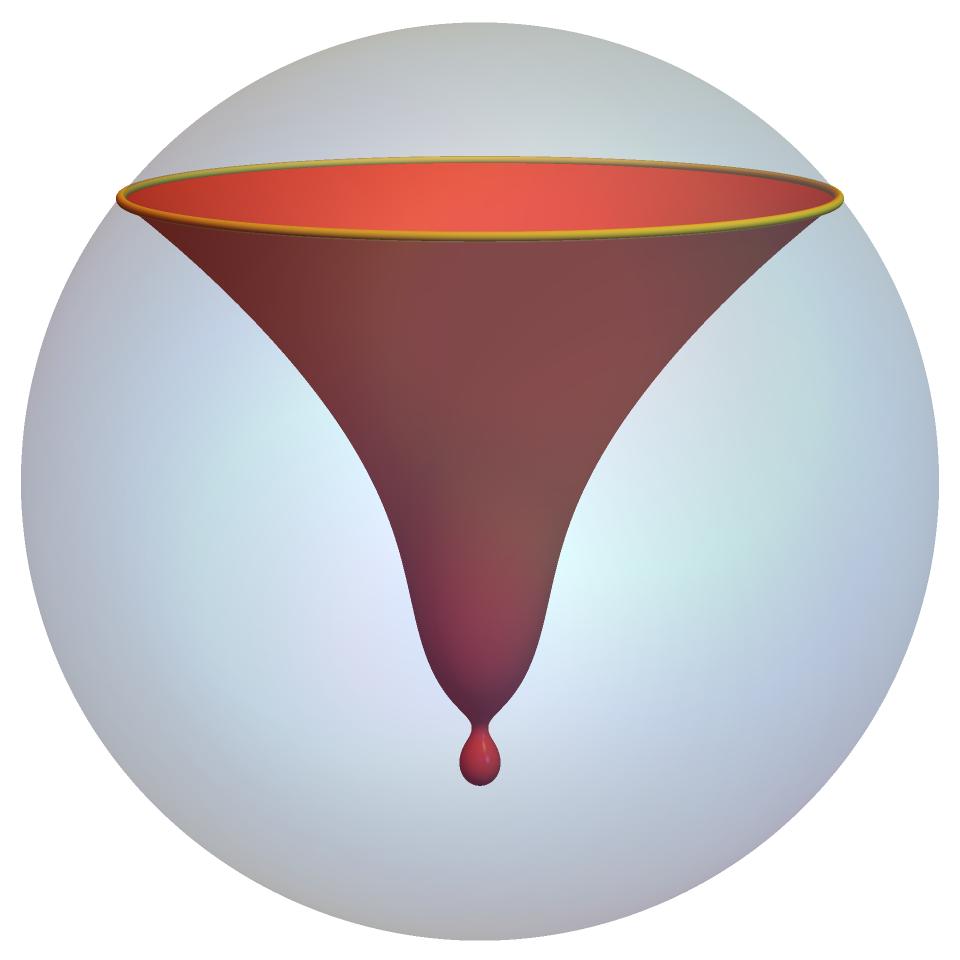}
		\end{subfigure}
		\,\,
		\begin{subfigure}[b]{0.4\linewidth}
			\includegraphics[width=\linewidth]{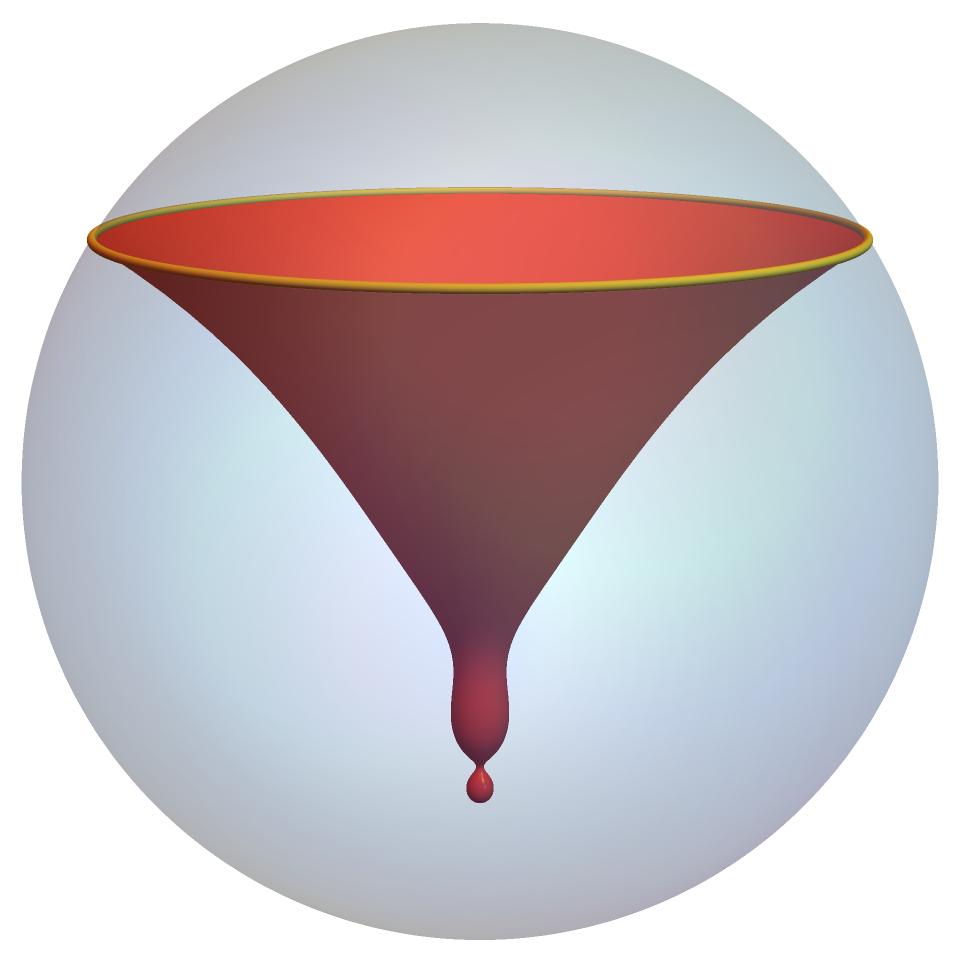}
		\end{subfigure}
		\,\,
		\begin{subfigure}[b]{0.4\linewidth}
			\includegraphics[width=\linewidth]{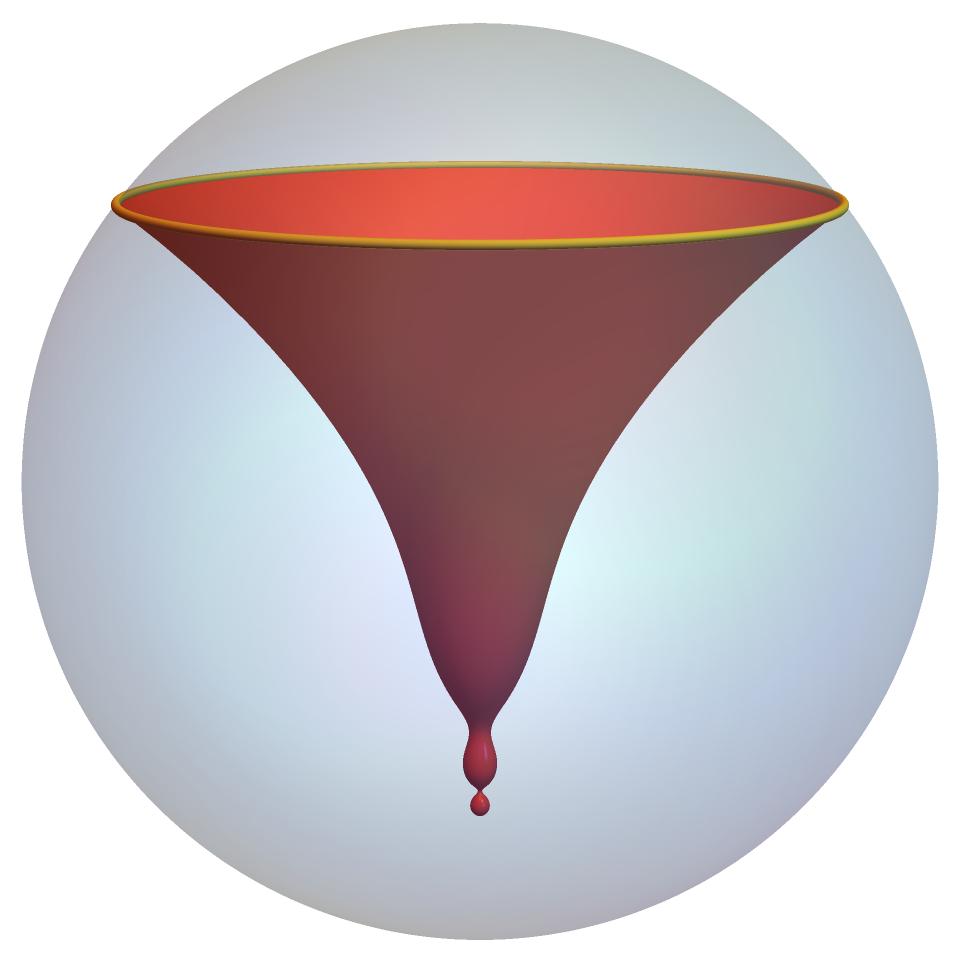}
	\end{subfigure}}
	\caption{{\small Six equilibria for the functional $-\mathcal{G}_R$ with $c_o=1$.}}
	\label{SurfacesHyp}
\end{figure}

In Figure \ref{SurfacesHyp} we illustrate (in the ball model for ${\bf H}^3$) six equilibrium surfaces for the functional $-\mathcal{G}_R$. To obtain these figures, we have numerically solved the system \eqref{ODE1}-\eqref{ODE3} with initial conditions \eqref{ic} for initial heights $z_o>0$ corresponding to the first six zeroes of the function illustrated on the left of Figure 3 of \cite{LPP}. We have then rotate the associated generating curves (from the top part until they meet $\{z=0\}$) around the $z$-axis and identify the upper half-space model with the ball model for ${\bf H}^3$, via the diffeomorphism
\begin{equation*}
	(x,y,z)\in{\bf R}^3_+\longmapsto \frac{1}{x^2+y^2+(z+1)^2}\left(2x,2y,x^2+y^2+z^2-1\right)\in{\bf B}^3.
\end{equation*}

We finish this section discussing that the minimum of $-\mathcal{G}_R$ among topological discs is not attained and obtaining a condition for stability of equilibria with respect to the functional\footnote{Observe that we are carrying the stability analysis for the functional $\mathcal{G}_R$, with the positive sign in front. Even though for equilibria we considered $-\mathcal{G}_R$ due to the relation with the Helfrich energy, our choice of sign here is more convenient for the stability analysis given the relationship between $\mathcal{G}_R$ and the hyperbolic area.} $\mathcal{G}_R$. 

We first notice that the minimum, among topological discs, of the functional $-\mathcal{G}_R$ for $c_o> 0$ is never attained. To the contrary, the minimum would be attained at an equilibrium surface $\Sigma$. From \eqref{inequality} and Theorem \ref{ELequations}, we have that equilibria for $-\mathcal{G}_R$ satisfy $-\mathcal{G}_R[\Sigma]=\mathcal{H}_{1,c_o,0}[\Sigma]$, which is non negative. However, the hemisphere $S_+^2$ of radius $1/c_o$ satisfies $-\mathcal{G}_R[S_+^2]=-2\pi<0$ (cf., \eqref{GRS}), contradicting our assumption that $\Sigma$ attains the minimum. In the case $c_o=0$, one can easily deform a surface $\Sigma$ while fixing a neighborhood of its boundary to produce a sequence of surfaces $\{\Sigma_i\}$ for which $-\mathcal{G}_R[\Sigma_i]\equiv-\mathcal{A}_R[\Sigma_i]\longrightarrow-\infty$.

We next study the second variation of the functional $\mathcal{G}_R$ to deduce a condition for stability. An equilibrium surface for the functional ${\mathcal G}_R$ will be called \emph{stable} if the second variation of the functional is non negative for all admissible variations. A necessary condition for stability will be that the second variation is non negative for all compactly supported normal variations $\delta X=\psi \nu$. 

Following \cite{PP3}, we define the surface invariant $\xi:=H+\nu_3/z$, so that the condition $\xi\equiv -c_o$ must hold for equilibrium. The linearization of this quantity defines a second order differential operator as the pointwise variation
$$P[\psi]:=2\,\delta\xi\,,$$
for normal variations $\delta X=\psi\nu$.

Following steps similar to those in Page 22 of \cite{PP4}, one obtains that for compactly supported normal  variations $\delta X=\psi\nu$,
\begin{equation}\label{SV}
	\delta^2 {\mathcal G}_R[\Sigma]=-\int_\Sigma \frac{\psi P[\psi]}{z^2}\:d\Sigma\:.
\end{equation}
In the following result we deduce that stable axially symmetric surfaces for $\mathcal{G}_R$ must be star shaped with respect to the origin of ${\bf R}^3$. In other words, for every point in the surface the line segment joining this point to the origin must be fully contained in the region enclosed by the surface and the plane $\{z=0\}$.

\begin{prop} 
	Let $\Sigma$ be an axially symmetric equilibrium surface for the functional ${\mathcal G}_R$ with $c_o>0$. If $\Sigma$ is stable, then $\Sigma$ is star shaped with respect to the origin.
\end{prop}
{\it Proof.\:} We will show that axially symmetric equilibria that are not star shaped with respect to the origin are unstable for $\mathcal{G}_R$. Assume that the equilibrium surface $\Sigma$ is not star shaped with respect to the origin. It then follows that the support function $q:=X\cdot\nu$ must vanish somewhere on the surface. If $q=0$ at some point in the surface, we let $z_1$ be the largest height of the circle at which this zero occurs and denote by $\overline{\Sigma}$ the part of $\Sigma$ above this circle. Note that $z_1>0$ holds since the position vector is parallel to the normal on the boundary circle located at $z=0$. Observe also that $q>0$ holds on $\overline{\Sigma}$ given our choice of unit normal.

As in Page 10 of \cite{PP3}, we compute the variation in the quantity $\xi=H+\nu_3/z$ with respect to the one parameter family of rescalings $X \longmapsto 
(1+\epsilon) X$. Under this change, we get $\xi\longmapsto (1+\epsilon)^{-1} \xi$. For this variation, $\delta X=q\nu +(X-q\nu)$ holds and it follows from the definition of $P$ given above that $P[q]=2c_o$ holds. 

We replace $q$ with the function $q_1$ which is equal to $q$ for $z\ge z_1$ and is identically zero for $z<z_1$. Therefore, for the variation $\delta X=q_1\nu$, the second variation formula \eqref{SV} gives
$$ \delta^2 {\mathcal G}_R[\Sigma]=-\int_\Sigma \frac{q_1 P[q_1]}{z^2}\,d\Sigma=-\int_{\overline{\Sigma}} \frac{q P[q]}{z^2}\,d\Sigma=-2c_o\int_{\overline{\Sigma}} \frac{q}{z^2}\,d\Sigma <0\:.$$
This shows that $\Sigma$ is unstable for $\mathcal{G}_R$. {\bf q.e.d.}
\\

We can observe that many  of the surfaces in Figure \ref{SurfacesHyp} (see Figure \ref{EHSurfaces} for their representation in ${\bf R}^3$) are not star shaped with respect to the origin so they are unstable for the functional ${\mathcal G}_R$. However, this does not imply anything about their stability for the Helfrich functional.

\section{Applications}

In this section we will apply our findings regarding equilibria for the functional $-\mathcal{G}_R$ to three different variational problems involving the Helfrich energy \eqref{H}.

\subsection{Closed Helfrich Surfaces.} Let $\Sigma$ be a closed surface of genus zero. The surface $\Sigma$ is a Helfrich topological sphere if and only if \eqref{ELH} holds on $\Sigma$. Observe that only the parameter $c_o\in{\bf R}$ affects the shape of closed equilibria for the Helfrich energy \eqref{H}. Indeed, the total Gaussian curvature term of $\mathcal{H}_{a,c_o,b}$ is a topological invariant by the Gauss-Bonnet Theorem and so the constant $b\in{\bf R}$ does not arise in the characterization of equilibria. In addition, $a>0$ is a multiplicative factor which does not play a role either.

In Theorem 4.1 of \cite{LPP}, axially symmetric Helfrich topological spheres with non-constant mean curvature were characterized as solutions of the reduced membrane equation \eqref{RME} which, in addition, satisfy the rescaling condition
\begin{equation}\label{rescaling}
	\int_\Sigma\left(H+c_o\right)d\Sigma=0\,.
\end{equation}

Disc type equilibria for $-\mathcal{G}_R$ can be smoothly reflected across $\{z=0\}$ to give closed genus zero surfaces. We next prove that these surfaces are Helfrich. 

\begin{prop}\label{CHS}
Let $\Sigma$ be an axially symmetric equilibrium disc type surface for the functional $-\mathcal{G}_R$ and denote by $\widetilde{\Sigma}$ its reflection across the plane $\{z=0\}$. Then, $\Sigma\cup\widetilde{\Sigma}$ is a closed genus zero axially symmetric Helfrich surface in ${\bf R}^3$.

Conversely, every axially symmetric Helfrich topological sphere which is symmetric with respect to the reflection through the plane $\{z=0\}$ arises in this way.
\end{prop}
{\it Proof.\:} Let $\Sigma$ be an axially symmetric disc type surface in equilibrium for $-\mathcal{G}_R$. From Theorem \ref{ELequations}, it satisfies the reduced membrane equation \eqref{RME}, and $\partial_nH=0$ holds at $\{z=0\}$.

We denote by $\widetilde{\Sigma}$ the reflection of $\Sigma$ across $\{z=0\}$ (which also satisfies \eqref{RME} since reflections across the horizontal plane $\{z=0\}$ preserve this equation) and construct the genus zero surface $\Sigma\cup\widetilde{\Sigma}$. Since $\partial_nH=0$ holds along $\partial\Sigma=\partial\widetilde{\Sigma}$, Corollary 3.12 of \cite{LPP} then shows that $\Sigma\cup\widetilde{\Sigma}$ is of class $\mathcal{C}^3$ at $z=0$. This regularity condition is equivalent to the surface being in equilibrium for the Helfrich energy (cf., Lemma 4.4 of \cite{LPP}). This proves the forward implication.

For the converse, let $\Sigma$ be an axially symmetric Helfrich topological sphere symmetric with respect to $\{z=0\}$ and denote by $\Sigma_+=\Sigma\cap\{z\geq 0\}$ and by $\Sigma_-=\Sigma\cap\{z\leq 0\}$. From Theorem 4.1 of \cite{LPP}, the reduced membrane equation \eqref{RME} holds on $\Sigma_+$ (and on $\Sigma_-$). If $c_o=0$, axially symmetric spheres satisfying \eqref{RME} are round. Hence, $\Sigma_+$ is a hemisphere and so an equilibrium for $-\mathcal{G}_R$. We next assume $c_o\neq 0$ holds. In this case, Theorem 4.1 of \cite{LPP} also shows that the rescaling condition \eqref{rescaling} is satisfied. Since $\Sigma_+$ and $\Sigma_-$ are reflections of each other, it follows that
\begin{eqnarray*}
	0&=&\int_\Sigma \left(H+c_o\right)d\Sigma=\int_{\Sigma_+}\left(H+c_o\right)d\Sigma+\int_{\Sigma_-}\left(H+c_o\right)d\Sigma\\
	&=&2\int_{\Sigma_+}\left(H+c_o\right)d\Sigma\,.
\end{eqnarray*}
From $c_o\neq 0$ and Lemma 4.5 of \cite{LPP} (see equation (22)), the mean value of $H+c_o$ vanishing on $\Sigma_+$ is equivalent to $\partial_n H=0$ at $z=0$. This, together with \eqref{RME}, shows that $\Sigma_+$ is in equilibrium for $-\mathcal{G}_R$. {\bf q.e.d.}
\\

Axially symmetric closed Helfrich surfaces of genus zero symmetric with respect to $\{z=0\}$ can be found in Figure 2 of \cite{LPP} (see also the last two pages of the same paper). The top half of these surfaces are equilibria for $-\mathcal{G}_R$. Indeed, the halves of the first six examples correspond with the surfaces of Figure \ref{SurfacesHyp} (after the identification of the ball and upper half-space models for ${\bf H}^3$).

We point out here that in \cite{LPP} numerical evidence was presented which suggests that, in fact, all axially symmetric Helfrich topological spheres are symmetric with respect to the plane $\{z=0\}$. The validity of this conjecture will allow us to remove this condition from the converse part of Proposition \ref{CHS}.

\subsection{Euler-Helfrich Surfaces.}

For a compact surface $\Sigma$ with smooth boundary, the \emph{Euler-Helfrich energy} functional is the Helfrich energy \eqref{H} on the interior of the surface coupled to the bending energy of the boundary. Thus, it is given by the total energy
\begin{equation}\label{E}
	E[\Sigma]:=\int_\Sigma\left(a\left[H+c_o\right]^2+bK\right)d\Sigma+\oint_{\partial\Sigma}\left(\alpha\kappa^2+\beta\right)ds\,,
\end{equation}
where $a>0$, $\alpha>0$, $\beta>0$ and $c_o$, $b$ are any real constants. Here, $\kappa$ denotes the (Frenet) curvature of the boundary $\partial\Sigma$.

The Euler-Lagrange equations characterizing equilibria for $E$ consist of \eqref{ELH} on the interior of $\Sigma$ together with the boundary conditions (for details, see Section 2 and Appendix A of \cite{PP1})
\begin{eqnarray}
	a(H+c_o)+b\kappa_n&=&0\,,\label{bc1}\\
	J'\cdot\nu-a\partial_nH+b\tau_g'&=&0\,,\label{bc2}\\
	J'\cdot n+a(H+c_o)^2+bK&=&0\,,\label{bc3}
\end{eqnarray}
where $J$ is the vector field along $\partial\Sigma$ given by
\begin{equation}\label{J}
	J:=2\alpha T''+\left(3\alpha\kappa^2-\beta\right)T\,.
\end{equation}
Here, $T$ is the unit tangent vector field to $\partial\Sigma$ and $\left(\,\right)'$ denotes the derivative with respect to the arc length parameter of the boundary.

We will next show that axially symmetric equilibria for $-\mathcal{G}_R$ are examples of Euler-Helfrich surfaces for suitable choices of the energy parameters.

\begin{prop}\label{EHS} Let $\Sigma$ be an axially symmetric equilibrium disc type surface for the functional $-\mathcal{G}_R$. Then, $\Sigma$ is an equilibrium surface for the Euler-Helfrich functional $E$, where the energy parameters are chosen so that the radius of the boundary circle is $r=\sqrt{\alpha/\beta}$ and
	$$2ac_o\sqrt{\alpha/\beta}=a+b\,,$$
holds.
\end{prop}
{\it Proof.\:} Let $\Sigma$ be an axially symmetric disc type surface in equilibrium for $-\mathcal{G}_R$. From Theorem \ref{ELequations}, the reduced membrane equation \eqref{RME} holds. In addition, the boundary circle is a geodesic ($\kappa_g\equiv 0$) along which $\partial_nH=0$ and $\tau_g\equiv 0$ holds (see Proposition \ref{propboundary}).

Since \eqref{RME} holds on $\Sigma$, the equation \eqref{ELH} is also satisfied. Moreover, the boundary conditions \eqref{bc1}-\eqref{bc3} simplify to
\begin{eqnarray*}
	a(H+c_o)+b\kappa_n&=&0\,,\\
	J'\cdot \nu&=&0\,,\\
	J'\cdot n+a(H+c_o)^2+bK&=&0\,.
\end{eqnarray*}
Observe that, as computed in Proposition \ref{propboundary} (see also the proof of Theorem \ref{axsym}),
\begin{equation}\label{Hz=0}
	H-\kappa_n-c_o=0\,,
\end{equation}
holds at $z=0$. Using this, we deduce that, along the boundary circle,
\begin{eqnarray*}
	a(H+c_o)^2+bK&=&a(H+c_o)^2+b\kappa_n(2H-\kappa_n)\\
	&=&a(H+c_o)^2+b\kappa_n(H+c_o)\\
	&=&(H+c_o)(a[H+c_o]+b\kappa_n)\\
	&=&0\,.
\end{eqnarray*}
The last equality follows from \eqref{bc1}. Therefore, the boundary conditions \eqref{bc2} and \eqref{bc3} further simplify to $J'\cdot\nu=J'\cdot n=0$. A simple computation from \eqref{J} then shows that the radius of the boundary circle must satisfy $r=\sqrt{\alpha/\beta}$. \emph{The vector field $J$ is the Noether current associated to translational invariance of the elastic energy of the boundary. The equations $J'\equiv 0$ characterize elastic curves. In the case of an elastic circle, the radius must satisfy $r=\sqrt{\alpha/\beta}$.}

Finally, we use \eqref{Hz=0} again to rewrite \eqref{bc1} as
\begin{eqnarray*}
	0&=&a(H+c_o)+b\kappa_n=a(\kappa_n+2c_o)+b\kappa_n\\
	&=&(a+b)\kappa_n+2ac_o=-(a+b)/r+2a c_o\,,
\end{eqnarray*}
since $\kappa_n=-1/r$ holds from the orthogonality of the intersection of $\Sigma$ with the plane $\{z=0\}$. This finishes the proof. {\bf q.e.d.}
\\

\begin{figure}[h!]
	\centering
	\includegraphics[height=7.2cm]{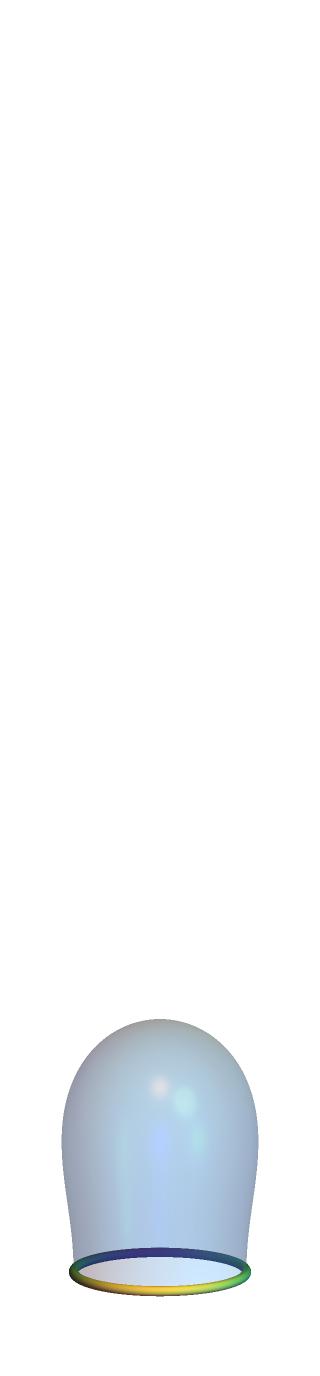}\quad
	\includegraphics[height=7.2cm]{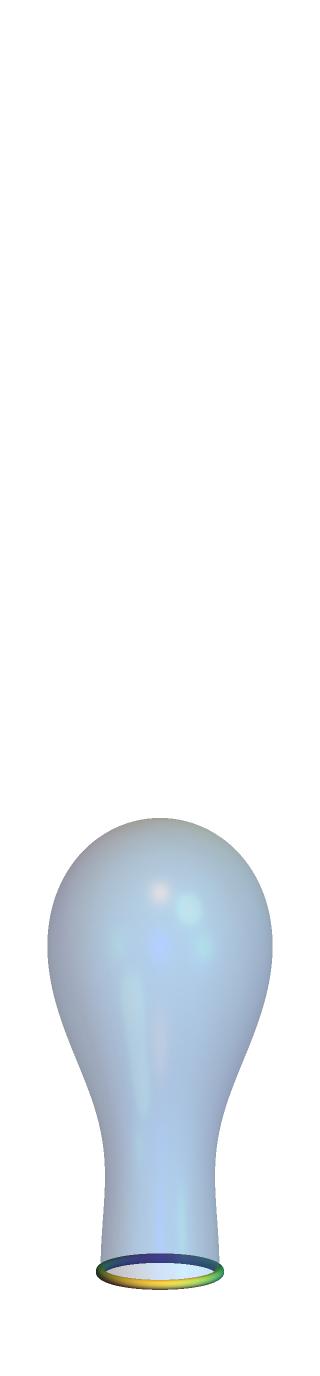}\quad
	\includegraphics[height=7.2cm]{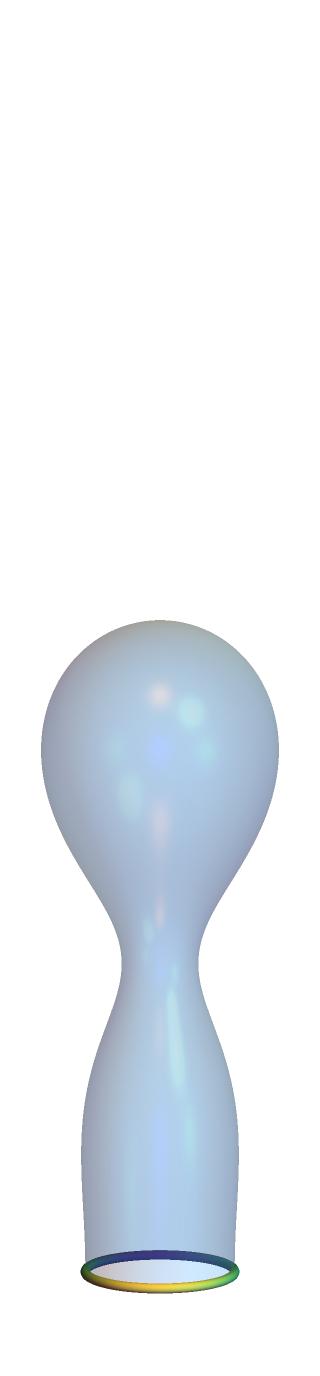}\quad
	\includegraphics[height=7.2cm]{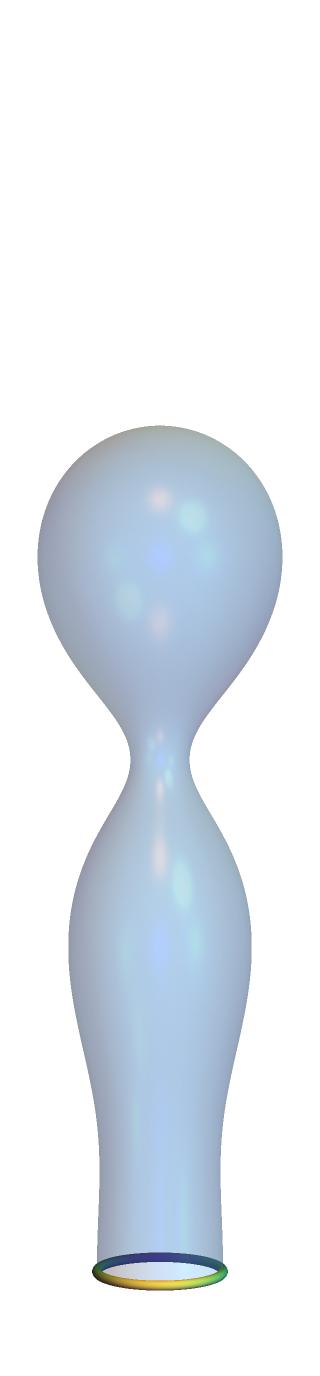}\quad
	\includegraphics[height=7.2cm]{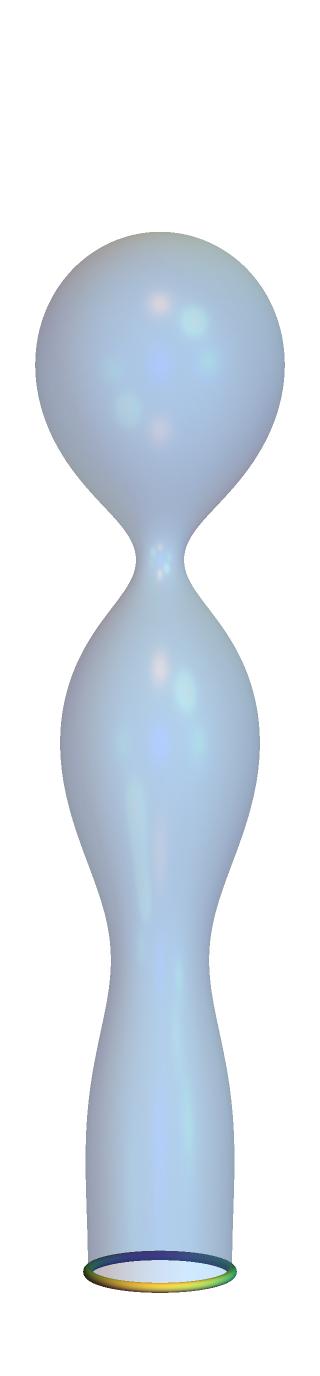}\quad
	\includegraphics[height=7.2cm]{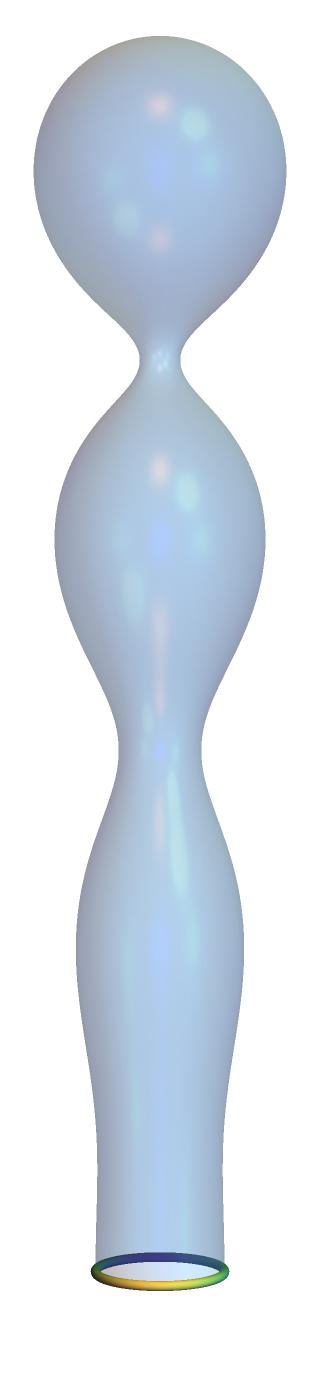}
	\caption{{\small The six equilibria for $-\mathcal{G}_R$ of Figure \ref{SurfacesHyp} illustrated in the upper half-space ${\bf R}^3_+$. These surfaces are in equilibrium for the Euler-Helfrich energy $E$ for suitable choices of the parameters.}}
	\label{EHSurfaces}
\end{figure}

In Figure \ref{EHSurfaces} we show the six equilibrium surfaces for the functional $-\mathcal{G}_R$ of Figure \ref{SurfacesHyp}, but this time represented in the upper half-space model for the hyperbolic space ${\bf H}^3$. When these surfaces are regarded as Euclidean surfaces in ${\bf R}^3_+$, they become equilibria for the Euler-Helfrich functional $E$ for a choice of parameters satisfying the relations of Proposition \ref{EHS}. 

\subsection{Helfrich Surfaces with Free Boundary.}

A Helfrich surface with free boundary is an equilibrium surface for the Helfrich energy \eqref{H} for \emph{admissible} variations. In the free boundary problem setting, admissible variations are those restricted to be tangent along the boundary to a fixed surface, known as the supporting surface, but are otherwise arbitrary. Here, we concentrate in the case that the supporting surface is the plane $\{z=0\}$.

Helfrich surfaces with free boundary in the plane $\{z=0\}$ are characterized\footnote{The Euler-Lagrange equations can be computed using standard techniques such as those of \cite{P}. More details will be presented in a future work. If preferred, for the moment the reader may just think of Helfrich surfaces with free boundary in $\{z=0\}$ as solutions of \eqref{ELH} satisfying \eqref{freebc1}-\eqref{freebc2}.} by \eqref{ELH} on the interior and
\begin{eqnarray}
	a(H+c_o)+b\kappa_n&=&0\,,\label{freebc1}\\
	\partial_nH&=&0\,,\label{freebc2}
\end{eqnarray}
along the boundary.

Axially symmetric equilibria for $-\mathcal{G}_R$ are Helfrich surfaces with free boundary in $\{z=0\}$, for suitable choices of the energy parameters for $\mathcal{H}_{a,c_o,b}$.

\begin{prop}\label{freeH} Let $\Sigma$ be an axially symmetric equilibrium disc type surface for the functional $-\mathcal{G}_R$. Then, $\Sigma$ is a Helfrich surface with free boundary in the plane $\{z=0\}$, where the energy parameters $a>0$ and $b$ are chosen so that the radius $r$ of the boundary circle satisfies
	$$2ac_o r=a+b\,.$$
\end{prop}
{\it Proof.\:} Let $\Sigma$ be an axially symmetric disc type surface in equilibrium for $-\mathcal{G}_R$. From Theorem \ref{ELequations}, the reduced membrane equation \eqref{RME} holds on $\Sigma$ and so does \eqref{ELH}. Moreover, $\partial_nH=0$ is satisfied along $\partial\Sigma$. Hence, \eqref{freebc2} automatically holds.

It only remains to rewrite \eqref{freebc1}. For this, we use \eqref{Hz=0} and the same argument as that of the last paragraph of the proof of Proposition \ref{EHS}. {\bf q.e.d.}
\\

\begin{figure}[h!]
	\centering
	\includegraphics[height=5cm]{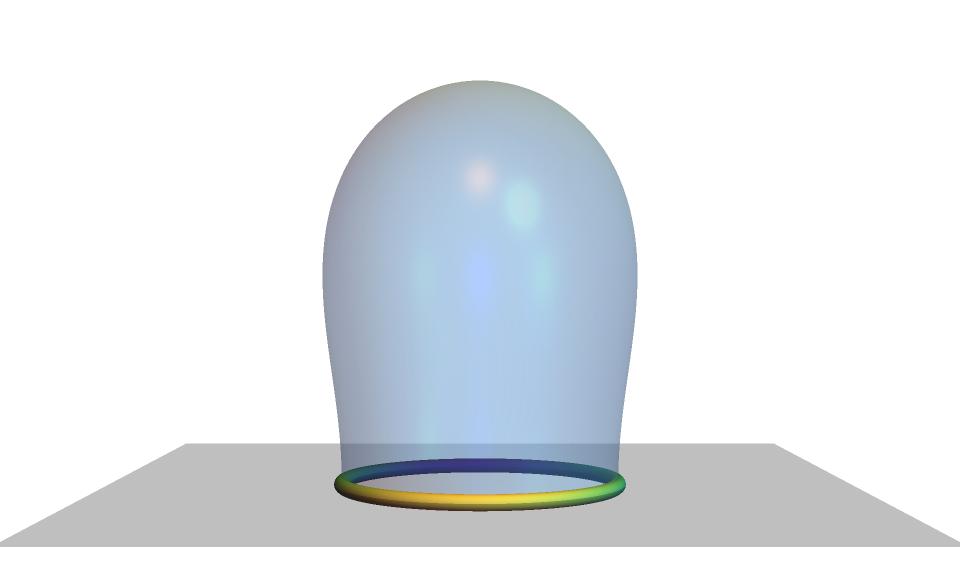}
	\caption{{\small A Helfrich surface with free boundary in the plane $\{z=0\}$. This surface is also the first equilibrium for the functional $-\mathcal{G}_R$ shown on Figures \ref{SurfacesHyp} and \ref{EHSurfaces}.}}
	\label{FBHSurfaces}
\end{figure}

To illustrate the result of Proposition \ref{freeH} we consider the first surface of Figure \ref{SurfacesHyp} (see its identification in ${\bf R}^3_+$ on the left of Figure \ref{EHSurfaces}). This surface is in equilibrium for $-\mathcal{G}_R$ and, when regarded as an Euclidean surface, it is also a Helfrich surface with free boundary in the plane $\{z=0\}$. We show this surface together with its supporting plane in Figure \ref{FBHSurfaces}.

\section*{Acknowledgments} 

The second author would like to thank Aaron J. Tyrrell for fruitful discussions regarding the renormalized area.





\bigskip

\begin{flushleft}
	Bennett P{\footnotesize ALMER}\\
	Department of Mathematics and Statistics,
	Idaho State University,
	Pocatello, ID 83209,
	U.S.A.\\
	E-mail: palmbenn@isu.edu
\end{flushleft}

\bigskip

\begin{flushleft}
	\'Alvaro P{\footnotesize \'AMPANO}\\
	Department of Mathematics and Statistics, Texas Tech University, Lubbock, TX
	79409, U.S.A.\\
	E-mail: alvaro.pampano@ttu.edu
\end{flushleft} 

\end{document}